\input mssymb
\magnification=\magstep1

\font\bbf=cmr10 scaled\magstep2
\def\a{\alpha}
\def\ad{{\rm ad}\,}
\def\ai{\a_i}
\def\aii{a_{ii}}

\def\b{\beta}
\def\be{$$}
\def\bm{\breve{m}}
\def\ca{\hat{\a}}
\def\cA{\breve{A}}
\def\cai{\hat{\a}_i}
\def\caii{\hat{a}_{ii}}
\def\caij{\hat{a}_{ij}}

\def\caji{\hat{a}_{ji}}
\def\cajj{\hat{a}_{jj}}

\def\ce{\hat{\e}}
\def\cft{conformal field theory}
\def\cfts{conformal field theories}
\def\cG{\breve{G}}
\def\cH{\hat{H}}
\def\char{{\rm ch}}
\def\chi{\Chi}
\def\Chi{\hat{h}_i}
\def\chj{\hat{h}_j}
\def\cht{\breve{{\rm ht}}}
\def\cI{\breve{I}}
\def\cIr{\breve{I}_{\rm r}}
\def\cl{\hat{\l}}
\def\cL{\breve{\L}}

\def\crho{\breve{\rho}}
\def\cW{\breve{W}}
\def\d{\delta}
\def\D{\Delta}
\def\D+{\Delta^{\!+}_{\phantom{X}}}
\def\diag {{\rm diag}}
\def\do{\dot{\omega}}
\def\e{\epsilon}

\def\eE{{\rm e}}
\def\g{\gamma}
\def\gl{{\rm gl}}
\def\go{\tilde\o}
\def\h{{\rm ht}}
\def\he{\hat{e}}
\def\hf{\hat{f}}

\def\imply{\Rightarrow}
\def\iN{\!\in\!}
\def\Ir{I_{\rm r}}
\def\kma{Kac-Moody algebra}
\def\l{\lambda}
\def\L{\Lambda}
\def\llb{\hbox{\bbf(}}
\def\lrb{\hbox{\bbf)}}
\def\m{\mu}
\def\mult{{\rm mult}}
\def\N{{\cal V}} 
\def\o{\omega}
\def\oaii{\sum_{l=0}^{N_i-1}\!a_{i,\do^li}}
  
\def\oaij{\sum_{l=0}^{N_j-1}\!a_{i,\do^lj}}
\def\oaji{\sum_{l=0}^{N_i-1}\!a_{j,\do^li}}
\def\oajj{\sum_{l=0}^{N_j-1}a_{j,\do^lj}}
  
\def\oai{\sum_{l=0}^{N_i-1}\alpha_{\do^li}}
\def\ohi{\sum_{l=0}^{N_i-1}h_{\do^li}}
\def\ohj{\sum_{l=0}^{N_j-1}h_{\do^lj}}
\def\or{\hat{r}}

\def\oW{\hat W}
\def\p{P_{\o}}
\def\pf{\hbox{\vbox{\hrule\hbox{\vrule\phantom{\vrule height 8pt width 8 pt
 depth 0pt}\vrule}\hrule}}}
\def\Ro{R^{\o}}
\def\sgma{{h^{}_{\tilde\o}}}
\def\tf{\tilde f}
\def\to{\rightarrow}
\font\eightrm=cmr8
\def\tom{\tilde\o}
\def\tr{{\rm tr}}

\def\z{\zeta}
\parindent 40pt
\hfuzz 2pt
\bigskip
\hfill DESY 96-098 

\hfill IHES/M/96/32

\hfill q-alg/9605046
\vskip2.8em
\centerline{\bbf Some automorphisms of Generalized Kac-Moody algebras}
\bigskip
\centerline{\bf J\"urgen Fuchs$^{\,1}$, Urmie Ray
\footnote{$^2$}{{\eightrm Supported by the EPSRC}}, Christoph 
Schweigert$^{\,3}$}
\bigskip
\centerline{{}$^1$\it{ DESY, 
Notkestra{\ss}e 85, D -- 22603 Hamburg}} \smallskip
\centerline{{}$^2$\it{ DPMMS, University of Cambridge, 16 Mill Lane, 
Cambridge CB2 1SB, U.K.}} \smallskip
\centerline{{}$^{2,3}$\it{ IHES, 35 Route de Chartres, F -- 91440 
Bures-sur-Yvette}}
\vskip2.8em

\noindent
{\bf 1. Introduction}
\bigskip\bigskip
In this paper we consider some algebraic structures associated to a class
of outer automorphisms of generalized Kac-Moody (GKM)
algebras. These structures have recently been
introduced in [2] for a smaller class of outer automorphisms in the case of 
ordinary Kac-Moody algebras with symmetrizable Cartan matrices. 

A GKM algebra $G=G(A)$ is essentially described by its Cartan matrix,\break
$A=(a_{ij})_{i,j\in I}$; the index set $I$ can be either a finite or a countably
infinite set. For any permutation $\do $ of the set $I$ which has finite 
order and leaves the Cartan matrix invariant, we find a 
family of outer automorphisms $\omega$ of the GKM algebra $G(A)$ which 
preserve the Cartan decomposition. 

Such an outer automorphism gives rise to a linear bijection $\tau_{\o}$ of 
$G$-modules, obeying the $\o$-twining property, i.e. if $V$ is a $G$-module, 
then
  $$\tau_\o (xv)=(\o^{-1}x)\tau_{\o}(v) $$
for all $x\iN G$ and all $v\iN V$. Thus in general $\tau_{\o}$ is not a 
homomorphism of $G$-modules, but some sort of ``twisted
homomorphism''. Furthermore $\tau_{\o}$ maps highest weight $G$-modules to
highest weight $G$-modules, though the image and pre-image are not always 
isomorphic.

In applications 
in conformal field theory, one is particularly
interested in those highest weight modules which are mapped to themselves.
The dual map $\o^*$ of the restriction of $\o$ to a Cartan
subalgebra $H$ of $G$ is a bijection of $H^*$, the dual of $H$.
For a highest weight $G$-module $V(\L)$ of highest weight $\L$ we have
$\tau_{\o}(V(\L))=V(\L)$ if and only if $\o^*(\L)=\L$.
A convenient tool to keep track of some properties of such a linear map is the
twining character of $V(\L)$, defined as in [2] as the formal sum
  $$(ch\, V)^{\o}=\sum_{\l\leq\L} m_{\l}^{\o}\,e(\l)\,,$$ where
  $$m_{\l}^{\o}=\cases{0\,,&if $\o^*(\l)\not=\l\,$;\cr 
  \noalign{\vskip2pt}
                   \tr({\tau_{\o}}|_{V_{\l}}^{})\,,&if $\o^*(\l)=\l\,$.\cr}$$

The main result of this paper is an explicit formula for the twining character
of Verma and irreducible highest weight $G$-modules. This formula shows that 
the twining characters can be described in terms of the characters
of highest weight modules of some other GKM algebra which depends on 
$G=G(A)$ and $\do $, the so-called orbit Lie algebra. In this paper we show 
that the `linking condition' that had to be imposed in [2]
is not needed, and that in particular this result applies to all Kac-Moody
algebras with symmetrizable Cartan matrices. 
In the case of affine Lie algebras,
this result has allowed for
the solution of two long-standing problems in \cft: the resolution
of field identification fixed points in diagonal coset \cfts\ (see [3]) and 
the resolution of fixed points in integer spin simple current modular 
invariants (see [4]).

In Section 2, after recalling the 
definition of a GKM algebra, we introduce the notion of an orbit Lie 
algebra and a twining character. In Section 3 we state and prove our main 
theorem, Theorem 3.1, which asserts that twining
characters are described by ordinary characters of the orbit Lie algebra, 
for a particular type of automorphisms of $G$ which just permute the
generators associated to simple roots.
As a by-product, we associate in Proposition 3.3 to the permutation $\do$
an interesting subgroup $\hat W$ of the Weyl group $W$ of a GKM algebra, 
which is again a Coxeter group. In Section 4 we extend the Theorem
to the whole class of outer automorphisms associated to a given finite order
permutation $\do$ of the index set $I$, leaving 
the Cartan matrix invariant.

Apart from the extension to arbitrary GKM algebras and to a larger class
of automorphisms, the present paper improves the treatment in [2] inasfar
as the description of $\hat W$ and the analysis of the cases with $\oaii\le 0$ 
are concerned. The corresponding statements, which previously had to be 
verified by detailed explicit calculations (see e.g.\ the appendix of [2]),
are now immediate consequences of our general results.
\bigskip \bigskip
\noindent
{\bf 2. Definitions and elementary properties}
\bigskip\bigskip
We first remind the reader of the definition of a Generalized Kac-Moody
(GKM) algebra and of some of its properties (see [1] or [5] for details).
All vector spaces considered are complex.  
Let $I$ be either a finite or a countably infinite set.  For
simplicity of notation, we identify $I$ with $\{1,2,...\,,n\}$ or $\Bbb Z_+$.
\smallskip\noindent
Let $A=(a_{ij})_{i,j\in I}$ be a matrix with real entries defined as follows:
\smallskip
  (i)$\,$ \ $a_{ij}\leq 0\,$ if $i\not=j$; \smallskip
  (ii)$\,$  ${{2a_{ij}}\over {a_{ii}}} \in\Bbb Z\,$ if   $a_{ii}>0$; \smallskip
  (iii)$\,$ if $a_{ij}=0$, then $a_{ji}=0$; \smallskip
  (iv)$\,$  there exists a diagonal matrix $D=\diag(\e_1,...,\e_n)$, with 
            $\e_i\iN\Bbb R$ and $\e_i>0$ 

            $\quad\,\,\ $ for all $i$, such that $DA$ is symmetric.
\medskip\noindent
A matrix satisfying condition (iv) is said to be symmetrizable.
\medskip
Let $H$ be an abelian Lie algebra of dimension greater or equal to $n$.
Let $h_1,...\,,h_n$ be linearly independent elements of $H$.
Define $\a_j$ in $H^*$, the dual of $H$, to be such that $\a_j(h_i)=a_{ij}$.  
\medskip
The GKM algebra $G=G(A)$ with Cartan matrix $A$ and Cartan subalgebra $H$ is a 
Lie algebra generated by $e_i,\, f_i\,$, $i\iN I$, and $H$,
with the following defining relations: \smallskip

  $[e_i, f_j]=\d_{ij}h_i \,,$ \smallskip
  $[h, e_i]=\a_{i}(h)e_i \,,$ \smallskip
  $[h, f_i]=-\a_{i}(h)f_i \,,$ \smallskip
  $(\ad e_i)^{1-2{a_{ij}/a_{ii}}}e_j=0=(\ad f_i)^{1-2{a_{ij}/ a_{ii}}} 
  f_j$\quad if $a_{ii}>0 \,,$ \smallskip
$[e_i, e_j]=0=[f_i, f_j]$\quad if $a_{ij}=0 \,.$
\medskip\noindent
{\bf Remarks.}  1. For any $n\times n$ diagonal matrix $D^{\prime}$ with
real positive entries, $G(D^{\prime}A)$ 
is a GKM algebra isomorphic to $G(A)$.  The resulting generators are
scalar multiples of the above ones, and so the simple roots $\a_j$ remain
unchanged. If we take $D^{\prime}$ to be 
$\diag(\e^{\prime}_1,\e^{\prime}_2,...\,,\e^{\prime}_n)$, where
  $$\e^{\prime}_i=\cases{2(a_{ii})^{-1},&if $a_{ii}>0\,$,\cr 
  \noalign{\vskip2pt}
              \ 1\,,&otherwise$\,$,\cr}$$
then all the diagonal positive entries of the matrix $D^{\prime}A$ are equal 
to 2, and we get a Cartan matrix as defined in [5].

\noindent
2. In [2], $h_i$, and $\a_j$ are defined to be so that $\a_j(h_i)=a_{ji}$.
In this paper we have taken the transpose in order to follow the 
convention in [5].
\medskip\noindent
The elements $h_1,h_2,...\,,h_n$ form a basis for $H\cap [G,G]$, and 
there is a subalgebra $C$ consisting of commuting derivations of $G$ such
that $H=H\cap [G,G]\oplus C$.  

\noindent
There exists a bilinear form $(.\,,.)$ on $H$ defined by:

  $$(h_i,h)=\e_i^{-1}\a_i(h)\,,\quad\hbox{and}\quad (h,h^{\prime})=0\,,\quad 
  h,h^{\prime}\in C\,,$$
where as above the matrix $\diag(\e_1,...\,,\e_n)A$ is symmetric.

\noindent
Therefore $(h_i,h_j)=\e_j^{-1}a_{ij}$.
This form extends uniquely to a bilinear, symmetric, invariant form on $G$,
whose kernel is contained in $H$. When the form is non-degenerate, it
induces a bilinear form on $H^*$, which we also denote by $(.\,,.)$.
In particular, $$(\a_i,\a_j)=\e_ia_{ij}\,.$$  Note that $H$ 
(and hence $G$) can always be
extended by adding outer derivations of $G$ having the $e_i$, $f_i$ as
eigenvectors so  as to make the form non-degenerate. 
When $(.\,,.)$ is non-degenerate on $H$, the Cartan decomposition holds
for the GKM algebra $G$:
  $$G= \llb\!\bigoplus_{\a\in\D+} G_{-\a}\lrb \oplus H\oplus
  \llb\!\bigoplus_{\a\in\D+} G_{\a}\lrb \,,$$
where $G_{\a}=\{x\iN G\mid [h,x]=\a(h)x,\, h\iN H\}$, and
$\D+$ is the set of positive roots of $G$ (i.e. $\a\iN H^*$ is a positive root
if $G_{\a}\not=0$, and $\a$ is a sum of simple roots).
\smallskip\noindent
If $a_{ii}>0$, the simple root $\a_i$ is called real. Set
  $$\Ir:=\{i\iN I\mid a_{ii}>0\}\,.$$
The Weyl group $W$ is generated by the reflections $r_i$, $i\iN\Ir$, acting on 
$H^*$. A root is said to be real if it is conjugate to a real simple root under 
the action of $W$, and imaginary otherwise. The group $W$ is a Coxeter group
(see Proposition 3.13 in [5]).  Recall that a Coxeter 
group is a group of the following type:
  $$<x_1,x_2,...\,,x_n\mid x_i^2=1;\, (x_ix_j)^{m_{ij}}=1\ (i,j=1,2,...\,,n,\, 
  i\neq j)>\,,$$ 
where the $m_{ij}$ are positive integers or $\infty$.
Set $m_{ij}$ to be the order of $(r_ir_j)$ for $i,j\iN\Ir$ ($i\neq j$).  
These orders are given by the following table, which we call $T$ (see [5]): 
  $$\vbox{\settabs 
  \+${2a_{ij}\over a_{ii}}{2a_{ji}\over a_{jj}}$\qquad & 
\quad $\geq 4$ \quad & \quad $\geq 4$ \quad & \quad $\geq 4$ \quad & 
\quad $\geq 4$ \quad & \quad $\geq 4$ \quad \cr 
  \+${2a_{ij}\over a_{ii}}{2a_{ji}\over a_{jj}}$&0&1&2&3&$\geq 4$\cr \smallskip 
  \+$m_{ij}$&2&3&4&6&$\infty$\cr}$$
\bigskip
In the rest of this paper $G=G(A)$ will denote a GKM algebra, with
non-degenerate bilinear form $(.\,,.)$.
Without loss of generality, we assume that the Cartan matrix $A$ is symmetric. 
\smallskip\noindent
Choose a bijection $\do\!:\, I\mapsto I$ of finite order which keeps 
the Cartan matrix fixed,\hfill\break
i.e.\ $a_{\do i, \do j}=a_{i,j}$ for all $i,j\iN I$.  
\smallskip\noindent
If the Dynkin diagram of $G$ is defined to be the Dynkin diagram of the GKM
subalgebra of $G$ generated by $e_i$, $f_i$ for all $i\in\Ir$,  
then $\do$ restricts to a bijection of the Dynkin diagram. 
(Note that the number of bonds linking node $i$ and $j$ is $\,{\rm max}\{{2
|a_{ij}|\over a_{ii}}\,,{2|a_{ji}|\over a_{jj}}\}$).
\smallskip\noindent
Let
$N$ be the order of $\do$ and $N_i$ the length of the $\do$-orbit of $i$ in $I$.
\smallskip\noindent
By the same arguments that were given in \S 3.2 of [2] for
Kac-Moody algebras, $\do$ induces an outer automorphism $\o$ of $G$ 
(the details of its action on the outer derivations in $H$ are given in [2]).  
In particular,
  $$\o e_i=e_{\do i}\,,\quad \o f_i=f_{\do i}\,,\quad \o h_i=h_{\do i}\,,
  \quad\hbox{for all}\ \ i\iN I\,.$$  
The automorphism $\o$ preserves the Cartan decomposition.
Let $\z\iN\Bbb C$ be a primitive $N$th root of unity. Then the
eigenvalues of the restriction $\o|_{H}$ of
$\o$ to $H$ are contained in $\{\z^l\mid l=0,1,...\,,N-1\}$.
Since $\o|_H$ has finite order, $H$ is the direct sum of its eigenspaces.
Let $H^l$ denote the eigenspace corresponding to eigenvalue $\z^l$.
\medskip\noindent
We choose a set of representatives from each $\do$-orbit:
  $$\hat I:=\{i\iN I\mid i\leq \do^li,\, \forall l\iN\Bbb Z\}\,.$$
Some of these orbits play a major role, so we also introduce the following 
subset of $\hat I$:
  $$\cI:=\{i\iN{\hat I}\mid\oaii\leq 0\imply \oaii=a_{ii}\} \,.$$
For $i\iN {\hat I}$ define 
  $$s_i:=\cases{{a_{ii}/\oaii}\,,&if $i\iN\cI$ and $a_{ii}\not=0\,,$\cr
  \noalign{\vskip 2pt}
  \ 1\,,&otherwise$\,$.\cr}$$
\smallskip\noindent
{\bf Remarks.}  1.  Suppose that $i\iN\cI$. 
The definition of $\cI$ implies that if $\a_i$ is imaginary, then 
$s_i=1$, and if  $\a_i$ is real, then we only have two possibilities: 
either $s_i=1$ and for all integers $1\leq l\leq N_i-1$, $a_{i,\do^li}=0$, or 
else $s_i=2$ and there is a unique integer $1\leq l\leq N_i-1$ such that 
$a_{i,\do^li}\not=0$. In the latter case, ${2a_{i,\do^li}\over a_{ii}}=-1$. So
we can deduce that $N_i$ is even, and $l={N_i\over 2}$.  
Hence when $s_i=1$, the orbit of $i$ in the Dynkin diagram of $G$ is totally
disconnected, i.e.\ of type $A_1\times\cdots\times A_1$ (where $A_1$ appears 
    $N_i$ times); and when $s_i=2$, 
the orbit of $i$ is of type $A_2\times\cdots
\times A_2$ (where $A_2$ appears $N_i\over 2$ times).

\noindent
2.  If $G$ is a Kac-Moody algebra and $\do$ fulfills the linking condition
of [2], then $\cI=\hat I$.  We will not need to impose this condition.
\smallskip
Define the matrix ${\hat A}=({\caij})_{{i,j}\in {\hat I}}$ to be as follows:
  $$\caij:=s_j\oaij.$$ 
\bigskip
\proclaim Lemma 2.1. The matrix $\hat A$ satisfies conditions {\rm(i)}, 
{\rm(ii)}, {\rm(iii)}, and {\rm(iv)} of the Cartan matrix of a GKM algebra.
\par\noindent
{\bf Proof.}
Suppose $i\not=j\in {\hat I}$.  Then $\caij\leq 0$ since for all integers $l$,
$\do^lj\not=i$ as $i$ and $j$ are not in the same $\do$-orbit.  
Suppose further that $\caii>0$. Then $\caii=\aii$, so that $2\caij\over\caii$ 
is an integer since $s_j$ is an integer for all $j\iN{\hat I}$.

\noindent
If $\caij=0$, then $0=\oaij={N_j\over N_i}\oaji$, so that $\caji=0$.
\hfill\break Let ${\hat D}=\diag(N_is_i)_{i\in{\hat I}}$.  Then 
straightforward calculations show that
${\hat D}{\hat A}$ is symmetric. This completes the proof.
\hfill\pf
\bigskip\noindent
Therefore  there is a GKM algebra, which we call $\hat G$, with Cartan matrix
$\hat A$, and such that the bilinear form induced by $\hat A$ is
non-degenerate on its Cartan subalgebra
$\cH$.  We let ${\he}_i$, and ${\hf}_i$ denote its other generators. Set
${\Chi}=[{\he}_i, {\hf}_i]\,,\ i\iN{\hat I}$.
\medskip\noindent
{\bf Remarks.} 1.  If $i\in\hat I-\cI$, then  ${a_{ii} / \oaii}$ is not
an integer when $a_{ii}\leq 0$, and it is non-positive when $a_{ii}>0$.
Therefore if $\hat I\not=\cI$, then the matrix with entries
$({a_{jj}\over\oajj})\sum_{l=0}^{N_j-1}
a_{i\do^lj}$ is not the Cartan matrix of a GKM algebra.

\smallskip\noindent 
2.  The elements of $G$ fixed by $\o$ form a GKM subalgebra of $G$
(see [1] for the proof).
This fixed point subalgebra has a GKM subalgebra whose Cartan matrix has
$(i,j)$-th entry equal to $\sum_{l=0}^{N_i-1}a_{\do^{l}i,j}$.  
However, $\hat G$ is not in general isomorphic to a subalgebra of $G$.
\medskip
We are now ready to define the {\it orbit Lie algebra\/} associated to the
automorphism $\o$. To do so, we have to use the subset $\cI$ rather than 
$\hat I$. 
\bigskip
\proclaim Definition 2.1.  
The orbit Lie algebra associated to the bijection $\do$ of the Cartan
matrix $A$, or equivalently to the automorphism $\o$ of $G$, is defined to
be the Lie subalgebra $\cG$ of $\hat G$ generated by ${\he}_i$, ${\hf}_i$ 
for $i\iN\cI$, and $\cH$.
\par\noindent
Lemma 2.1 implies that $\cG$ is a GKM algebra with Cartan matrix
  $${\cA}=({\caij})_{{i,j}\in {\cI}}\,.$$
\medskip\noindent
{\bf Remarks.} 1.  The GKM subalgebra $\cG$ of $\hat G$ is also 
in general not isomorphic to a subalgebra of $G$.

\noindent
2.  The set $\cI$ may be empty, in which case $\cG=0$.

\noindent
3.  It can be shown that if $G$ is of finite type, then so
is the orbit Lie algebra $\cG$; and if $G$ is of affine (resp.\ indefinite)
type, then $\cG$ is either trivial, or also of affine (resp.\ indefinite) 
type (see [2]).
\medskip
We next define a linear map $P_{\o}\!:\, H^0\cap [G,G]\rightarrow\cH$ as
follows:
  $$\p\llb\ohi\lrb=N_i\Chi\,.$$ 
\bigskip
\proclaim Lemma 2.2.  For all $h,h^{\prime}\iN [G,G]\cap H^0$,
$\,(h,h^{\prime})=(\p(h),\p(h^{\prime}))$.
\par\noindent
{\bf Proof.}  Let $i$ be in $\hat I$. Since $N_i=N_{\do^li}$ for all 
integers $l$,
  $$\eqalign{(\ohi,\ohj)&=N_i\oaij =N_i s_j^{-1} \caij\cr
  &=(N_i\Chi,N_j\chj)=(P_{\o}\ohi,P_{\o}\ohj)\,.\cr}$$
The result follows by linearity.
\hfill\pf
\smallskip\noindent
Provided we choose $\cH$ to have the right number of outer derivations,
this map can be extended to the outer derivations contained in $H^0$ so as to
give an isomorphism $H^0\rightarrow\cH$,
in such a way that Lemma 2.2 holds for all $h,h^{\prime}$ in
$H^0$ (the proof is the same as in \S 3.3 of [2], where this is shown for 
Kac-Moody algebras). For simplicity of notation, we also call this isomorphism 
$\p$.
\smallskip\noindent
The automorphism $\o$ induces a dual map $\o^*$ on $H^*$, namely:
  $$(\o^*\b)(h)=\b(\o h)\,,\quad\hbox{for}\quad \b\iN H^*,\,\, h\iN H.$$
In particular,
  $$\o^*(\a_i)=\a_{\do^{-1}i}\,,$$
since $(\o^*(\a_i))(h) = \a_i (\o(h)) =(h_i,\o(h))=
(h_{\do^{-1}i},h)=\a_{\do^{-1}i}(h)$ for all $h\iN H$.

\noindent
This bijection has the same eigenvalues as $\o|_H$, and so $(H^*)^l$ will
denote the eigenspace corresponding to eigenvalue $\z^l$. 
\smallskip\noindent
Since its restriction to $H^0$ is 
non-degenerate, the bilinear form gives rise to a bijection between
$(H^0)^*$ and $(H^*)^0$. Hence $\p$ induces a dual map $\p^*\!:\,\cH^*
\rightarrow (H^*)^0$. By definition $\l\iN (H^*)^0$ if and only if
$\o^*(\l)=\l$. Such weights will be called {\it symmetric weights}. Set
  $$\b_i:=\oai\,,\quad\hbox{for each}\quad i\in I.$$
In particular the following holds:
\bigskip
\proclaim Lemma 2.3. \item{{\rm(i)}} \ For all $i\iN{\hat I}$, 
$\p^*(\cai)=s_i\b_i\,$; \ and
\vskip2pt
\item{{\rm(ii)}} \ for all $\l,\mu\iN (H^*)^0,\,$ $(\l,\mu) = (\p^{*-1}(\l),
\p^{*-1}(\mu))\,.$
\par\noindent
{\bf Proof.}  For all $h\iN H^0$, Lemma 2.2 implies that 
  $$\ca_i(\p(h))=N_is_i\,(\chi,\p(h))
  =s_i\,(\ohi,h)\,,$$ 
so that (i) holds.

\noindent
(ii) follows directly from Lemma 2.2 and the definition of $\p^*$.
\hfill\pf
\bigskip\noindent
We next define the {\it twining character\/} of a highest weight $G$-module.
We first need to associate a representation $R^{\o}$ to a given representation
$R$ of the GKM algebra $G$.
\bigskip
\proclaim Definition 2.2.  Let $V$ be a $G$-module and $R\!:G\to\gl(V)$
the corresponding representation.  Define $R^{\o}$ to be the representation
$G\to\gl(V)$ such that $R^{\o}(x)=R(\o(x))$ for all $x\iN G$.
\par\smallskip\noindent
Let $R_{\L}\!:\,G\to \gl(V(\L))$ be a 
highest weight $G$-representation of highest weight $\L$ in $H^*$. 
Then $(R_{\L})^{\o}$ is a highest weight representation of highest weight
$\o^*(\L)$, $R_{\o^*(\L)}\!:\,G\to\gl(V(\o^*(\L)))$, since
$\o$ preserves the Cartan decomposition.  
Thus the automorphism $\o$ induces a bijection of $G$-modules
$\tau_{\o}\!:\, V({\L})\to V({\o^*(\L)})$ which satisfies the
{\it{$\o$-twining property}}, i.e.
  $$\tau_\o (R_{\L}(x)v)=R_{\o^*(\L)}(\o^{-1}x)\tau_{\o}(v)
  \quad{\rm for\ all}\ v\iN V(\L),\,\, x\iN G\,.$$
\medskip\noindent
{\bf Remarks.} When $\o^*(\L)\not=\L$, the representations $R_{\L}$ and 
$(R_{\L})^{\o}$ are not isomorphic. The bijection $\tau_{\o}$ is a linear map, 
but not in general an isomorphism of $G$-modules, even if $\o^*(\L)=\L$.
\medskip
   Denote by $V_{\l}:=\{v\iN V\mid hv=\l(h)v,\,h\iN H\}$ the weight space of 
$V$ of weight $\l$.  The bijection
$\tau_{\o}$ maps $V_{\l}$ onto $V_{\o^*(\l)}$ 
In particular if $R$ corresponds to the Verma (resp.\
irreducible highest weight) module $M({\L})$ (resp. $L({\L}$)), then 
$\Ro$ corresponds to the Verma (resp.\
irreducible highest weight) module $M({\o^*(\L)})$ (resp. $L({\o^*(\L)}$)).  
\smallskip\noindent
In this paper, we study the case when $\L$ is a symmetric weight (i.e. 
$\o^*(\L)=\L$), so
that $\tau_{\o}$ maps $M({\L})$ (respectively $L({\L})$) to itself. In the 
physics literature, such weights are called {\it fixed points\/} (see [6]).

The ordinary character $\char\,V$ of a highest weight $G$-module $V$
is the formal sum
  $$\char\,V:=\sum_{\l} (\dim V_{\l})\,e(\l)\,.$$
Replacing the formal exponential by the
exponential function, this gives a complex valued function
  $$\char\,V(h)=\sum_{\l}(\dim V_{\l})\,\eE^{\l(h)}=\tr_V^{} \eE^h\,,$$ 
defined on the set $Y(V)$ of elements $h\iN H$ such that the series converges 
absolutely.

Let $\rho\iN H^*$ be a Weyl vector for $G$, 
i.e. $(\rho,\a_i)={1\over 2}(\a_i,\a_i)$
for all $i\iN I$.  Such a vector exists since by assumption $(.\,,.)$ is
non-degenerate on $H$.
For $w\iN W$, let $\e(w)=(-1)^l$, where
$l$ is the minimal number of simple reflections $r_i$ needed to write $w$.
Let $S_{\L}=e(\L+\rho)\sum_{\b}\e(\b)e(-\b)$, where $\e(\b)=(-1)^m$ if $\b$ 
is a sum of $m$ distinct pairwise orthogonal imaginary simple roots,
orthogonal to $\L$, and $\e(\b)=0$ otherwise.
Borcherds showed that if $\L\iN H^*$, $(\L,\a_i)\geq 0$ for
all $i\iN I$, and $2(\L,\a_i)/(\a_i,\a_i)\iN\Bbb Z$ for all real simple
roots $\a_i$ of $G$, then the irreducible module $L(\L)$ of highest weight 
$\L$ has character
  $$\char\, L(\L)=\sum_{w\in W} \e(w)w(S_{\L})/e(\rho)\prod_{\a\in\Delta^+}
  (1-e(-\a))^{\mult\a}\,.$$
(For details of the proof of the character formula, see [1] or [5].)
\bigskip
\proclaim Definition 2.3. Let $\L\in H^*$ be a symmetric weight.  We 
define the twining character for 
the highest weight representation $R_{\L}$ on $V({\L})$ 
to be the following complex valued function defined on $Y(V)$:
$(\char\,V(\L))^{\o}(h) = \tr_{V}^{}\tau_{\o}\eE^{R_{\L}(h)}$.
\par\noindent
Note that since the twining character is bounded by the ordinary character, 
it is absolutely convergent on $Y(V)$.
Equivalently the twining character is the formal sum
  $$(\char\, V(\L))^{\o}=\sum_{\l\leq\L} m_{\l}^{\o}\,e(\l)\,,$$
where 
  $$m_{\l}^{\o}=\cases{0\,,&if $\o^*(\l)\not=\l\,$;\cr
  \noalign{\vskip 2pt}
                   \tr({\tau_{\o}}|_{V_{\l}}^{})\,,&if $\o^*(\l)=\l\,$.\cr}$$
Let $\N_{\L}$ (resp. $\Psi_{\L}$) denote the ordinary character of
the Verma (resp.\ irreducible) $G$-module of highest weight $\L$, and
$\breve{\N}_{\cL}$ (resp. $\breve{\Psi}_{\cL}$) denote the ordinary character of
the Verma (resp. irreducible) $\cG$-module of highest weight $\cL$.  
\medskip\noindent
{\bf Remarks.} A weight $\L$ in $H^*$ is said to be integrable if
$(\L,\a_i)\geq 0$ for all $i\iN I$ and $2(\L,\a_i)\over (\a_i,\ai)$
is an integer for all real simple roots $\a_i$.
A $G$-module $V$ is called integrable
if $f_i$ and $e_i$ act locally nilpotently for all $i\iN I$ such that
$a_{ii}>0$. An irreducible $G$-module of highest weight $\L$ is unitarizable
if and only if $\L$ is integrable.    
If all the simple roots of $G$ are real, i.e.\ $G$ is a Kac-Moody
algebra, then an irreducible $G$-module of
highest weight $\L$ is integrable if and only if it is unitarizable
(see \S 3, 10, and 11 in [5] for details).
\bigskip
\noindent
{\bf 3.  Twining characters and orbit Lie algebras}
\bigskip\noindent
We  now state the main result of this paper.
\bigskip
\proclaim Theorem 3.1.  Let $\L\in{H^*}$ be a symmetric weight, i.e.
$\o^*(\L)=\L$.
The twining character of the Verma 
$G$-module of highest weight $\L$ coincides with the ordinary character
of the Verma $\cG$-module of highest weight $\p^{*-1}(\L)$:
  $$\p^{*-1}(\N_{\L})^{\o}={\breve{\N}}_{\p^{*-1}(\L)}.$$
If, moreover, $\L$ is integrable, then 
the twining character of the irreducible 
$G$-module of highest weight $\L$ coincides with the ordinary character
of the irreducible $\cG$-module of integrable highest weight $\p^{*-1}(\L)$:
  $$\p^{*-1}(\Psi_{\L})^{\o}=\breve{\Psi}_{\p^{*-1}(\L)}.$$
\par\bigskip\noindent
In order to prove this Theorem, we first need a few more results.
Any Weyl vector $\rho$ of $G$ satisfies $(\o^*(\rho),\a_i)=(\rho,\a_i)$
   for all $i\iN I$, and hence we can choose 
$\rho$ to be a symmetric weight.
\bigskip
\proclaim Lemma 3.2.  The weight $\crho=P^{*-1}(\rho)$ is a Weyl vector in 
$\cH^*$ for $\cG$.
\par\noindent
{\bf Proof.}  Lemma 2.3 implies that if $i\iN\cI$ and $a_{ii}\not=0$, then 
  $$(\crho,\cai)=s_i(\rho,\b_i) ={1\over 2}s_iN_ia_{ii}
  ={1\over 2}s_i^2N_i\oaii ={1\over 2}(\ca_i,\cai)\,.$$  
The penultimate equality follows from the definition of $s_i$.

\noindent
If $i\iN\cI$ and $a_{ii}=0$, then $a_{i,\do^li}=0$ for all integers $l$,
so that $(\crho,\cai)=0={1\over 2}(\cai,\cai)$.
\hfill\pf
\medskip\noindent
{\bf Remarks.}  1. The proof of the previous Lemma shows that in order for
$\crho$ to be a Weyl vector for $\cG$, we need to scale the numbers
$\oaij$ by $s_j$ to define the Cartan matrix $\cA$.  

\noindent
2.  For $i\in{\hat I}-\cI$, one has
$(\crho,\cai)\not={1\over 2}(\ca_i,\cai)$, so that when ${\hat I}\not=\cI$, 
$\crho$ is {\it not\/} a Weyl vector for the bigger GKM algebra $\hat G$.  
It is not possible to scale the Cartan matrix $\hat A$ in such a way
that on the one hand, the resultant matrix remains the Cartan matrix of a 
GKM algebra, and on the other hand, 
the resultant vector $\crho$ is a Weyl vector for the corresponding GKM algebra.
\medskip\noindent
Let $W$ be the Weyl group of $G$, and $\cW$ the Weyl group of ${\cG}$.  
Note that $\caii>0$ for $i\iN{\hat I}$ implies that $i\iN\cI$. Since ${\hat G}$ 
and $\cG$ have the same Cartan subalgebra, $\cW$ is therefore also the Weyl
group of $\hat G$. Let
$$\cIr :=\{i\iN{\hat I}\mid \caii>0\}\,,$$ and for $i\iN\cIr$ let $\or_i$
denote the reflections corresponding to the
simple real roots of $\hat G$ (or equivalently $\cG$). Define
  $$\oW:=\{w\iN W\mid w\o^*=\o^*w\}$$ 
to be the set of all elements in the Weyl group $W$ of $G$ commuting 
with the bijection $\o^*$ of $H^*$. 
This is a subgroup of $W$. 

If $s_i=2$, then the orbit of $i$ in the Dynkin diagram of $G$ is the product 
of ${N_i\over 2}$ copies of the Dynkin diagram of $A_2$. Then 
$\{\do^li,\do^{l+{N_i\over 2}}i\}$ are the connected components of the orbit 
of $i$. For each $i\iN\cIr$ (i.e. $\caii>0$), define 
  $$w_i:=\cases{r_i r_{\do i} \cdots r_{\do^{N_i-1} i}\,,& if $s_i=1\,$;\cr 
\noalign{\vskip2pt}
    \prod_{l=1}^{{N_i/2}} r_{\do^li} r_{\do^{l+{N_i/2}}i}r_{\do^li}\,,&
if $s_i=2\,$.\cr}$$ 
As
in \S 5.1 of [2], it can be shown that the elements $w_i$ are in 
$\oW$, and that for symmetric weights $\l\iN H^*$, 
  $$w_i(\l)=\l-{2s_i\,(\l,\a_{i})\over (\a_i,\a_i)}\oai \,.\eqno(1)$$
In fact the elements $w_i$ generate the group $\oW$:
\bigskip
\proclaim Proposition 3.3.  $\,\oW =\,\,<w_i\mid i\iN\cIr>\,$.
\par\noindent
{\bf Proof.}  Set ${\tilde W}:=\,\,<w_i\mid i\iN\cIr>\,$. Let $\L\iN H^*$ 
be an integrable symmetric weight such that $(\L,\a_i)>0$ for all
$i\iN I$.  Such weights exist since $(.\,,.)$ is non-degenerate on $H^*$.  
Let $\l\leq\L$ be a symmetric weight.  
\smallskip
\noindent
We claim that if $i\in\Ir$ and $\b_i=\oai$ satisfies
$(\b_i,\b_i)\leq 0$, then $(\l,\a_i)\geq 0$.

\noindent
Since both $\L$ and $\l$ are symmetric, $\L-\l=\sum_{i\in{\hat I}}k_i\b_i$, 
where $k_i\geq 0$ for each $i\iN{\hat I}$. Since $(\b_i,\b_i)=N_i(\b_i,\a_i)$ 
and $(\b_j,\a_i)\leq 0$ for all $j\neq i$, our claim follows.
\smallskip\noindent
Since $\L$ is integrable, $w(\l)\leq\L$ for all $w$ in the Weyl group $W$
(see \S 3 in [5]).  
Let $w\in{\tilde W}$ be such that the height of $\L-w(\l)$ is minimal, i.e.\
$\h(\L-w(\l))\leq\h(\L-w^{\prime}(\l))$ for all $w^{\prime}\iN{\tilde W}$. 
\smallskip 
\noindent
We claim that $w(\l)$ is in the positive Weyl chamber, i.e.\
for all $i\iN I$ such that $a_{ii}>0$, $(w(\l),\a_i)\geq 0$.

\noindent
Assume this is false. Since $w$ commutes with $\o$, $w(\l)$ is symmetric.
Hence the above argument implies that for all $j\in\hat I_{\rm r}-\cIr$, 
$(w(\l),\a_j)\geq 0$.
Thus there is some $i\iN\cIr$ such that $(w(\l),\a_i)<0$. {}From $(1)$ we get
$w_iw(\l)=w(\l)-s_i{2(w(\l),\a_{i})\over (\a_i,\a_i)}\b_i$, so that 
$\h (\L-w_iw(\l))<\h (\L-w(\l))$, contradicting the definition of $w$.
\smallskip
\noindent
Let $w^{\prime}$ be in $\oW$.  Then ${w^\prime}(\L)$ is a symmetric weight.  So 
from the above, there is some $\tilde w$ in $\tilde W$ such that 
${\tilde w}w^{\prime}(\L)$ is in the positive Weyl chamber. Since the 
$W$-orbit of $\L$ intersects the positive Weyl chamber at a unique
point, we can deduce that ${\tilde w}w^{\prime}(\L)=\L$. Furthermore by 
definition $(\L,\a_i)\not=0$ for all $i\iN I$.
Hence ${\tilde w}w^{\prime}=1$ (see Proposition 3.12 in [5]),
  so that $w'\iN\tilde W$, and hence $\tilde W=\oW$.
\hfill\pf
\bigskip\noindent
The next result shows that $\oW$ is a Coxeter group.
\bigskip
\proclaim Corollary 3.4.  The subgroup $\oW$ of $W$ is isomorphic to 
the Weyl group $\cW$ of $\cG$.
\par\noindent
{\bf Proof.} We first show that for all $i,j\iN\cIr$, $(w_iw_j)^{\bm_{ij}}=1$,
where the exponents $\bm_{ij}$ are given by table $T$ (changing $a_{ij}$,
$a_{ji}$, $a_{ii}$, and $a_{jj}$ to $\caij$, $\caji$, $\caii$, and $\cajj$
respectively, in the table). 
{}From Lemma 2.3 and $\hat D=\diag(N_is_i)$, for all $i\iN\cIr$ and all 
$\cl\iN (\cH)^*$ we have 
  $${(\cl,\cai)\over (\cai,\cai)}= {(\p^*(\cl),\b_i)\over N_ia_{ii}}
  ={N_i\,(\p^*(\cl),\a_i)\over N_ia_{ii}}$$
since $\p^*(\cl)$ is symmetric and $\p^*(\cai)=s_i\b_i$. Therefore 
$\p^*(\or_i(\cl))=w_i(\p^*(\cl))$ follows by comparison with $(1)$. Now 
$\cW$ is the Coxeter group characterized by $(\or_i\or_j)^{\bm_{ij}}=1$.
So by induction on the number of generators 
$\or_i$ and $\or_j$ in
the expression $(\or_i\or_j)^{\bm_{ij}}$, we can deduce from what precedes that
  $$\p^*(\cl)=(w_iw_j)^{\bm_{ij}}(\p^*(\cl))\,.$$
Let $w:=(w_iw_j)^{\bm_{ij}}$.
Since $\p^*$ is a bijection between $\cH^*$ and $(H^*)^0$,
$w(\l)=\l$ for all symmetric weights $\l$ in $H^*$.
In particular $w(\rho)=\rho$ as $\rho$ is assumed to be symmetric.
The definition of $\rho$ implies that for all $i\iN\Ir$ (i.e.\ such that
$a_{ii}>0$), $(\rho,\a_i)>0$.  Hence the proof of Proposition 3.12 
in [5] tells us that $w=1$.
\smallskip\noindent
We may therefore define a map $\,\Theta\!:\, \cW\to\oW$ such that
  $$\Theta(\or_i)=w_{i}\,,$$  
which extends in a natural way to $\cW$.  
The above reasoning shows that $\Theta$ is well defined and a group 
homomorphism. Proposition 3.3 implies that $\Theta$ is surjective.

\noindent
It only remains to show that $\Theta$ is injective.
{}From the preceding calculations we can also deduce that
$$\p^*(\or(\cl))=\Theta(\or)(\p^*(\cl))\,,$$
for all elements $\or\iN\cW$.  So again the bijectivity of $\p^*$
implies that if $\Theta(\or)=1$, then $\or=1$.
Thus $\Theta$ is a group isomorphism.
\hfill\pf
\bigskip
As the next two results show, with 
respect to the twining character, the subgroup $\oW$ plays the role
that the Weyl group plays with respect to the ordinary character:
\bigskip
\proclaim Proposition 3.5.  If $V$ is an integrable 
highest weight $G$-module with highest weight $\L$, then
$\,w((\char\, V)^{\o})=(\char\, V)^{\o}$ for all $w\iN\oW$.
\par\noindent
{\bf Proof.}  Let $R$ be the representation: $G\to\gl(V)$, 
$\l$ a symmetric weight of $V$, and $V_{\l}$ the corresponding
weight space in $V$. Then $w_i(\l)$ is symmetric, and is
a weight of $V$ since $\L$ is integrable.
Set $x_l^R:=(\exp f_{\do^li})(\exp -e_{\do^li})(\exp f_{\do^li})$.  Define
 $$X_i:=\cases{x_1^Rx_2^R\cdots x_{N_i-1}^R,& if $s_i=1\,$;\cr
\noalign{\vskip 3pt}
 \prod_{l=1}^{{N_i/2}}x_l^Rx_{l+{N_i/2}}^Rx_l^R\,,&if $s_i=2\,$.\cr}$$

\noindent
Lemma 3.8 in [5] implies that $X_i(V_{\l})=V_{w_i(\l)}$. Since
$\o$ extends uniquely to an automorphism of the universal
enveloping algebra $U(G)$ of $G$, the definition 
of $\tau_{\o}$ implies that $$\tau_{\o}(X_iv)=(\o^{-1}X_i)\tau_{\o}(v)\eqno(2)$$
for all $v\iN V$.  Since $w_i$ commutes with $\o^*$ and the Coxeter
relations hold for $\oW$,\break 
$\o^{-1}X_i=X_i$. Hence
we can deduce from $(2)$ that the trace of $\tau_{\o}$ on $V_{\l}$ equals the
trace of $\tau_{\o}$ on $V_{w_i(\l)}$. Therefore $w_i((\char\, V)^{\o})
=(\char\, V)^{\o}$, and the result follows from Proposition 3.3.
\hfill\pf
\bigskip\noindent
For $w\iN\oW$, let $\hat l(w)$ be the minimal number of generators $w_i$ 
needed to write $w$. Define 
  $$\ce(w):= (-1)^{\hat l(w)}.$$
Let $\l\leq\L$ be symmetric weights in $H^*$.  Consider the Verma module
$M(\L)$.  
Taking a basis of the universal algebra $U(G)$ of $G$ as given by the PBW 
theorem, we can deduce that
the trace of $\tau_{\o}$ on $M(\L)_{\l}$ only depends on the action
of $\o$ on $U(G)$.  Therefore
the expression $e(-\L-\rho)\N_{\L}^{\o}$ is independent of $\L$. Set
  $$\N^{\o}:=e(-\L-\rho)\,\N_{\L}^{\o}.$$
\bigskip
\proclaim Lemma 3.6.  For all $w\iN\oW$, $\,w(\N^{\o})=\ce(w)\N^{\o}$.
\par\noindent
{\bf Proof.}  
Consider the Verma module $M(0)$ with highest weight $0$.
Let $i\iN\cI$, and  
  $$\Delta_i:=\{-\b\iN\Delta^+\mid \exists\, j\not=\do^l i\,\,\,\forall l\iN\Bbb Z
  \,,\ \hbox{such that}\ \a_j\leq\b\}\,. $$
The ordinary character of $M(0)$ is 
  $\N_0=\prod_{\a\in\Delta^+}(1-e(-\a))^{-\mult\a}\,.$
Therefore all weights $\m\iN H^*$ such that $\m\leq 0$ are weights of $M(0)$.
Let $\l\leq 0$ be a symmetric weight.  We can write
$\l=\sum_{\b\in\Delta_i}\b-n\oai$ for some 
non-negative integer $n$.
\smallskip\noindent
When $s_i=1$, then for all integers $l,l^{\prime}$,
$\a_{\do^li}+\a_{\do^{l^{\prime}}i}$ is not a root. 
So we can order the positive roots of $G$ in the following way:
$\a_i$, $\a_{\do i},\ldots,\a_{\do^{N_i-1}i},\g_1,\g_2,\ldots\,$, with
$\g_p\iN\Delta_i$.

\noindent
When $s_i=2$, then for $l,l^{\prime}\iN\Bbb Z$,
$\a_{\do^li}+\a_{\do^{l^{\prime}}i}$ is a root if and
only if $l^{\prime}\equiv l+{N_i\over 2}\ ({\rm mod}\, {N_i})$; and for all
$l,l^{\prime},l^{\prime\prime}\iN\Bbb Z$,
$\a_{\do^li}+\a_{\do^{l^{\prime}}i}+\a_{\do^{l^{\prime\prime}}i}$ is
not a root. In this case, we can order the positive roots of $G$ in the 
following way: $\a_i^{}$, $\a_{\do i}^{},\ldots,\a_{\do^{N_i-1}i},
\a_i^{}+\a_{\do^{N_i/2}i},\ldots, \a_{\do^{{N_i/2}-1}i}+\a_{\do^{N_i-1}i},
\g_1,\g_2,\ldots\,$, with $\g_p\iN\Delta_i$.
\smallskip\noindent
By definition, $\tau_{\o}(xv)=\o^{-1}(x)\tau_{\o}(v)$ for all
$x\iN U(G)$ and all $v\iN M(0)$.
So choosing a basis of $M(0)_{\l}$ given by the PBW theorem, which respects
the above ordering of roots, we can deduce that the only basis vectors
of $M(0)_{\l}$ contributing to the trace of $\tau_{\o}$ are as follows:
\smallskip\noindent
For the case $s_i=1$, these vectors are
  $$f_i^mf_{\do i}^m\cdots f_{\do^{N_i-1}i}^mv_q^{(m)}\,,$$ 
where $0\leq m\leq n$ and
the vectors $v_q^{(m)}$ form a basis of the weight space $M(0)_{\l+m\b_i}$;

\noindent
and for the case $s_i=2$, if $\tf_i:=[f_i,f_{\do^{N_i/2}i}]$, these vectors are
  $$f_i^{m_0}f_{\do i}^{m_1}\cdots f_{\do^{N_i-1}i}^{m_{N_i-1}}
  \tf_i^{n_0}\tf_{\do i}^{n_1}\cdots \tf_{\do^{N_i/2-1}i}^{n_{N_i/2-1}}
  v_q^{(m_k,n_j)}\,,$$
where $m_k=m_{k+N_i/2}$, $\,m_k+n_k=m_j+n_j\leq n$ for all
$0\leq j,k\leq N_i/2-1$,
and the vectors $v_q^{(m_k,n_j)}$ form a basis of 
the weight space $M(0)_{\l+(m_0+n_0)\b_i}$.
\smallskip\noindent
Commutator terms that arise when reshuffling the products of the
generators of the corresponding root spaces to the form given by the
chosen basis can never give rise to a non-zero contribution to the
trace of $\tau_{\o}$ in 
  $M(0)_{\l}$.
Therefore if $m_{\l}^{\o}:=\tr (\tau_{\o})|_{M(0)_{\l}}$, then
summing over all the symmetric weight spaces of $M(0)$, we can deduce that
  $$\N_0^{\o} = \llb\sum_{\l}m^{\o}_{\l}e(\l)\lrb \,(1+m_{-\b_i}^{\o}e(-\b_i)
  +m_{-2\b_i}^{\o}e(-2\b_i) + \ldots\,)\,,\eqno(3)$$
where the first sum is taken over all sums of roots in $\Delta_i$.

\noindent{}From the above we can also deduce that for $s_i=1$, 
  $$\tr(\tau_{\o})_{M(0)_{-n\b_i}} =1$$ 
since all the simple roots $\a_{\do^li}$ are pairwise orthogonal;
and for $s_i=2$,
  $$\tr(\tau_{\o})_{M(0)_{-n\b_i}}=
  \sum_{\scriptstyle0\le n_j\le n\atop\scriptstyle0\le j\le N_i/2-1}
  (-1)^{\sum_{j=0}^{N_i/2-1}n_j}= \llb\sum_{k=0}^n (-1)^k\lrb^{N_i/2-1},$$
so that 
  $$\tr(\tau_{\o})_{M(0)_{-n\b_i}}=\cases{1\,,&if $n\equiv 0\pmod 2\,,$\cr
  \noalign{\vskip2pt}                                               
                                               0\,,& otherwise$\,$.\cr}$$ 
Substituting the value of $\tr(\tau_{\o})_{M(0)_{-n\b_i}}$ in $(3)$, we get
  $$\N_0^{\o}= \llb\sum_{\l\in{\Delta}_i}m^{\o}_{\l}e(\l)\lrb
  \,(1-e(-s_i\b_i))^{-1} \,.$$
Since $\rho$ is symmetric by assumption, and $(\rho,\a_i)={1\over 2}(\a_i,\a_i)$
for all $i\iN I$, $(1)$ gives: $w_i(\rho)= \rho -s_i\b_i$.
Now $r_i$ permutes the set of all negative roots distinct from $-\a_i$.  
Therefore $w_i$ permutes the elements of ${\Delta}_i$ and 
multiplies $e(-\rho)\,(1-e(-s_i\b_i))^{-1}$ by $-1$.
Hence the assertion of the Lemma follows for the generator $w_i$. This 
argument works for each $i\iN\cIr$, and so Proposition 3.3
implies the Lemma for all $w\iN\hat W$.
\hfill\pf
\bigskip\noindent
Set
  $${B^{\o}_{\!\L}}:=\{\l\iN H^*\mid \l\leq\L,\, \bigl\vert\l+\rho\bigr\vert=
  \bigl\vert\L+\rho\bigr\vert,\, \o^*(\l)=\l\}.\,$$
Arguments similar to those used in [2] or in \S 9 of [5] imply that 
for any symmetric weight $\L$ in $H^*$,
  $$\N_{\L}^{\o}=\sum_{\l\in {B^{\o}_{\!\L}}}
  c_{\L\l}^{}\,\Psi_{\l}^{\o},$$ where $c_{\L\l}\in\Bbb C$,
and $c_{\L\L}=1$.  Since $B^{\o}_{\!\L}$ is a discrete set, by inverting
the upper triangular matrix $(c_{\L\l})_{\l\in B^{\o}_{\!\L}}$, we get
  $$\Psi_{\L}^{\o}=\sum_{\l\in {B^{\o}_{\!\L}}}
  c_{\l}^{}\,\N_{\l}^{\o},\eqno(4)$$ where $c_{\L}=1$.
\bigskip
\proclaim Lemma 3.7.  Suppose that $\L\iN H^*$ is a symmetric integrable 
weight. Then for each scalar $c_{\l}$ in equation $(4)$, there
is some $w\iN\hat W$ such that $c_{\l}=\ce(w)c_{w(\L-\a+\rho)-\rho}$,
where $\a$ is a symmetric sum of 
distinct pairwise orthogonal imaginary simple roots, all orthogonal to $\L$.
\par\noindent
{\bf Proof.} {}From $(4)$ we get
  $$\Psi_{\L}^{\o}=\N^{\o}\sum_{\l\in {B^{\o}_{\!\L}}}
  c_{\l}^{}\,e(\l+\rho)\,,\eqno(5)$$ 
where $\N^{\o}$ is independent of $\l$. Given $\l\iN B^{\o}_{\!\L}$,
let $w\iN\hat W$ be such that the height of $\L+\rho-w(\l+\rho)$ is minimal.
Let $\mu:=w(\l+\rho)-\rho$. The proof of Proposition 3.3 shows that 
$(\l+\rho,\a_i)\geq 0$ for all real simple roots $\a_i$.  
Then $\mu=\L-\sum_{i\in I}k_i\a_i$, where the $k_i$
are non-negative integers. Furthermore $\vert\mu+\rho\vert^2=
\vert\L+\rho\vert^2$ implies that 
  $$\sum_{i\in I} k_i(\L,\a_i)+\sum_{i\in I}k_i(\mu+2\rho,\a_i)=0\,.$$ 
So as in the 
proof of Theorem 11.13.3 in [5] it follows that if $k_i\not=0$ then $\a_i$
is imaginary and $(\a_i,\L)=0$; that $(\a_i,\a_j)=0$ if 
$k_j$ is also non-zero for $j\not=i$; and that $(\a_i,\a_i)=0$ if $k_i\geq 2$. 
Since $\Psi_{\L}^{\o}$ is bounded by the ordinary character
$\Psi_{\L}$, terms such as $e(\L-\sum_{i\in I}k_i\a_i)$ do not occur in
$\Psi_{\L}^{\o}$ as all the roots $\a_i$ are orthogonal to $\L$ (see \S 11 of 
[5]). The ordinary character of the Verma module $M(\L)$ equals
$e(\L)\prod_{\a\in{\Delta}^+}(1-e(-\a))^{-\mult\a}$, so that
$\N^{\o}$ is bounded by $e(-\rho)\prod_{\a\in{\Delta}^+}(1-e(-\a))^{-\mult\a}$.
Hence we can deduce that if $k_i\not=0$, then $k_i=1$.
The Lemma now follows from Proposition 3.5 and Lemma 3.6.
\hfill\pf
\bigskip
We next determine the values of the scalars $c_{\l}$ in
$(5)$ when $\L-\l$ is a sum of
imaginary simple roots, pairwise orthogonal, and all orthogonal to $\L$.
Note that since $\L$ and $\l$ are both symmetric, so is $\L-\l$.  Hence
$\l=\L-\sum_{i\in\cI}\b_i$, where the sum can be taken to be over $\cI$ rather
than the larger $\hat I$, as $(\a_i,\a_{\do^li})=0$ for all integers 
$1\leq l\leq N_i-1$ implies that $i$ is $\do$-conjugate to an integer in $\cI$.

\smallskip\noindent
We first need another definition. If $\gamma=\sum_{i\in\cI}k_i\b_i$, define
  $$\cht(\g):=\sum_{i\in\cI}k_i \,.$$
\bigskip
\proclaim Lemma 3.8.  Let $\L$ be a symmetric integrable weight in 
$H^*$, and $\l$ be a symmetric
element in ${B^{\o}_{\!\L}}$.  If $\L-\l$ is a sum of distinct, pairwise
orthogonal, imaginary simple roots, all orthogonal to $\L$, then
the coefficient $c_{\l}$ in {\rm(4)} is
  $$c_{\l}=(-1)^{\cht(\L-\l)}.$$
\par\noindent
{\bf Proof.}  We prove this Lemma by induction on $\cht(\L-\l)$.
We know from $(4)$ that $c_{\L}=1$. Set $C_{\L}$ to be the set of all symmetric
weights $\mu$ in $B^{\o}_{\!\L}$ such that $\L-\mu$ is the sum 
of distinct, pairwise orthogonal, imaginary simple roots, and let
  $$\N_{0}^{\o}=\sum m^{\o}_{\mu}\,e(\mu)\,.$$  
Since the ordinary character of the Verma module $M(0)$ of $G$ is
  $$\N_{0} = \prod_{\a\in{\Delta}^+}(1-e(-\a))^{-\mult\a}\,,$$ 
for all $\mu$ in $C_0$ the dimension of the weight space $M(0)_{\mu}$ is 1, 
as $-\mu$ is a sum of orthogonal simple roots. So the definition of 
$\tau_{\o}$ gives $m^{\o}_{\mu}=1$ for all $\mu\iN C_0$. Let $\l\iN C_{\L}$, 
and so $\L-\l=\sum_{s=1}^r\b_{i_s}$, where the $\b_{i_s}$ are all distinct.  
The coefficient of $e(\l)$ on the right hand side of $(4)$ is
  $$\sum_{s=0}^r\sum_{\{j_1,...,j_s\}\in T_s}c_{\L-\b_{i_{j_1}}-...
  -\b_{i_{j_s}}}\,,$$
where $T_s$ is the set of all subsets of $\{1,2,...\,,r\}$ of order $s$.
Since $e(\l)$ does not appear on the left hand side of $(4)$, this
coefficient equals 0.
Assume now that the result holds for all weights $\mu$ in $C_{\L}$ such that
$\cht(\L-\mu)<\cht(\L-\l)$.  It follows by induction that
  $$c_{\l}+(-1)^{\cht(\L-\l)}\sum_{s=1}^r\pmatrix{r\cr s\cr}(-1)^{s}=0\,,$$
which gives the desired answer for $c_{\l}$.
\hfill\pf
\bigskip\noindent
Combining the results of Lemma 3.7 and Lemma 3.8
we have therefore proved that when $\L$ is a symmetric integrable weight, then
  $$\Psi_{\L}^{\o}=\N^{\o}\sum_{w\in\oW}\ce(w)\,w(S_{\L}^{\o})\,,$$
where
  $$S_{\L}^{\o}=e(\L+\rho)\sum\ce(\b)\,e(-\b)\,,$$ 
and $\,\ce(\b)=(-1)^{\cht(\b)}$ if $\b$ is the symmetric sum of pairwise
orthogonal imaginary simple roots, all orthogonal to $\L$, and
$\ce(\b)=0$ otherwise. Also of course, for the trivial module
$\Psi_0^{\o}=1$; we can therefore deduce that
  $$\N^{\o}= \llb\!\sum_{w\in\oW}\ce(w)w(S_{0}^{\o})\lrb^{-1}\,.$$
Substituting this result back into the formula for $\Psi_{\L}^{\o}$, we obtain
for any integrable symmetric weight $\L$ in $H^*$,
  $$\Psi_{\L}^{\o}={\sum_{w\in\oW}\ce(w)\,w(S_{\L}^{\o})\over
  \sum_{w\in\oW}\ce(w)\,w(S_{0}^{\o})}\,,$$
and for any symmetric weight $\L$ in $H^*$
  $$\N_{\L}^{\o}=e(\L)\, \llb\!\sum_{w\in\oW}\ce(w)w(S_{0}^{\o})\lrb^{-1}\,.$$
We can now complete the Proof of Theorem 3.1.
\bigskip\noindent
{\bf Proof of Theorem 3.1.}
Let $i\iN{\hat I}$ and $\a_i$ be imaginary, then 
$(\a_i,\o^{*l}(\a_i))=0$ for all integers $1\leq l\leq N_i-1$ 
if and only if $i\iN\cI$ and $\cai$ is an imaginary root of $\cG$.

\noindent
When $\L$ is an integrable symmetric weight in $H^*$, it follows
from Lemma 2.3 that $\p^{*-1}(\L)$ is an integrable weight
in $\cH^*$ for the GKM algebra $\cG$ (note that
it is not integrable for the bigger algebra $\hat G$).  
Furthermore Corollary 3.4 implies that 
the minimal number of generators $w_i$ of $w\iN\oW$ equals the number of
generators $\or_i$ of $\Theta^{-1}(w)$ in $\cW$.  Therefore the ordinary
character of the irreducible $\cG$-module of highest weight $\p^{*-1}(\L)$
equals $\p^{*-1}\Psi_{\L}^{\o}$, when $\L$ is integrable; and 
for any symmetric weight $\L$, the character of the Verma $\cG$-module
of highest weight $\p^{*-1}(\L)$ equals $\p^{*-1}\N_{\L}^{\o}$.
This completes the proof of Theorem 3.1.
\hfill\pf
\bigskip\noindent
It follows that if $V$ is a highest weight $G$-module with symmetric highest
weight $\L$, and
  $$(\char\, V)^{\o}=\sum_{\l\leq\L} m^{\o}_{\l}\,e(\l)\,,$$ 
then $m^{\o}_{\l}\not=0$ implies that
$\L-\l=\sum_{i\in\cI} k_i\b_i$, where for all $i\iN\cI$, $k_i$ is a 
non-negative integer.  Note that the sum may be taken to be over $\cI$,
and not the larger $\hat I$.
The denominator formula for $\cG$ immediately gives the following Corollary.
\bigskip
\proclaim Corollary 3.9.  Let ${\breve\Delta}^{\!+}_{\phantom{X}}$ denote the
set of positive roots of $\cG$. If $\L$ is a symmetric weight in $H^*$, then
$\,\N_{\L}^{\o}=e(\L)\,\prod_{\ca\in{{\breve\Delta}^{\!+}_{\phantom{X}}}}
(1-e(\p^*(-\ca)))^{-\mult\ca}$, where $\,\mult\,\ca=\dim\cG_{\a}$.
\par\bigskip\noindent
{\bf Remarks.} 1.  When $\hat G\neq \cG$,
the twining characters for highest weight $G$-modules
coincide with ordinary characters of $\cG$ and not of $\hat G$.
This is due to the fact that when ${\hat I}\not=\cI$,
$\p^{*-1}(\rho)$ is not always a Weyl vector for $\hat G$.

\noindent
2.  If $\cI=\emptyset$, the above results implies that 
for all symmetric weights $\L$ in $H^*$, $\N_{\L}^{\o}=e(\L)$,
and $\Psi_{\L}^{\o}=e(\L)$ if $\L$ is also integrable.  In this case $\cG=0$, 
and $\p^{*-1}(\L)=0$, so that $\p^{*-1}(\N_{\L}^{\o})=e(0)$ and
$\p^{*-1}\Psi_{\L}^{\o}=e(0)$. This is Theorem 2 in [2]. 

\noindent
3.  Since for fixed Cartan matrix,
the ordinary character and the twining character do not depend on 
the size of the Cartan subalgebra, Theorem 3.1 holds for any GKM algebras 
$G=G(A)$ and $\cG=G(\cA)$ as long as the duals of their Cartan subalgebras 
are large enough for all roots to be either positive or negative, for the 
multiplicities of the simple roots to be finite, and for the existence of 
Weyl vectors.  In particular the bilinear forms induced by $A$ and $\cA$ 
need not be non-degenerate, and the multiplicities of the
simple roots may be greater than 1.
\bigskip\noindent
{\bf 4.  A larger class of outer automorphisms.}
\bigskip\noindent
As before the Cartan matrix $A$ is symmetric, $G=G(A)$ denotes
a GKM algebra such that the bilinear form on $G$ induced by $A$ is
non-degenerate; and the bijection $\do$ of the set $I$ preserves the
Cartan matrix $A$.  Also, $\o$ denotes the outer automorphism of $G$ defined
in section 2. To the bijection $\do$ we can in fact not only associate 
the automorphism $\o$, but a whole family of outer automorphisms of $G$. More 
precisely, there exist automorphisms $\tilde\o$ of $G$ such that
  $$\go e_i=\xi_i e_{\do i}\,,\quad\go f_i=\xi^{\prime}_i f_{\do i}\,,$$
where $\xi_i$, and $\xi^{\prime}_i$ are in $\Bbb C^*$.
(The proof of the existence of $\o$ in [2] applies to $\go$ as well).
It follows immediately that
  $$\go=\phi\o\,,$$ 
where $\phi$ is an inner automorphism of $G$ such that
  $$\phi e_i=\xi_i e_{i}\,,\quad\phi f_i=\xi^{\prime}_i f_{i}\,.$$
We will refer to the automorphisms $\o$ and $\go$ as diagram automorphisms
and generalized diagram automorphisms, respectively.
\bigskip\noindent
\proclaim Lemma 4.1.  With the above notation,
$\xi^{\prime}_i=\xi^{-1}_i$ for all $i\iN I$ for which there exists $j\iN I$
such that $a_{ij}\not=0$.
\par\noindent
{\bf Proof.}  Since $\phi$ is a Lie algebra homomorphism,
$\phi (h_i)=\xi_i\xi^{\prime}_ih_i$. Now on the one hand, for all $j\iN I$,
$\,\phi([h_i,e_j])=\phi(a_{ij}e_j)=a_{ij}\phi(e_j)$. On the other hand,
$[\phi(h_i),\phi(e_j)]=\xi_i\xi^{\prime}_ia_{ij}\phi(e_j)$.  
Hence $a_{ij}\phi(e_j)=\xi_i\xi^{\prime}_ia_{ij}\phi(e_j)$, and the result 
follows.
\hfill\pf
\bigskip
We now assume that for all $i\iN I$, $\xi^{\prime}_i=\xi^{-1}_i$.
Therefore $\o$ and $\tilde\o$ are equal on the Cartan subalgebra $H$, so that
$\o^*={\tilde\o}^*$ (i.e. the dual of ${\tilde\o}|_H$ on $H^*$).
As in the previous section, $\tilde\o$ induces a bijection $\tau_{\tilde\o}$
of $G$-modules: $\tau_{\tilde\o}\!:\, V_{\L}\to V_{\o^*(\L)}$, which 
satisfies the {\it{$\tilde\o$-twining property}}, i.e.
  $$\tau_{\tilde\o} (R_{\L}(x)v)=R_{{\tilde\o}^*(\L)}(\go^{-1}x)
  \tau_{\tilde\o}(v) $$
for all $x\iN G$ and all $v\iN V_\L$.
When $\o^*(\L)=\L$, we define the twining character of a 
$G$-module $V$ with respect to $\tilde\o$ to be
  $$(\char\, V)^{\tilde\o}(h) = \tr_{V}^{}\tau_{\tilde\o}\eE^{R_{\L}(h)}.$$
We next show that $(\char\, V)^{\tilde\o}$ can be easily expressed in terms of
$(\char\, V)^{\o}$.

Note that unlike in section 3, we now require 
the bilinear form on $H$ to be non-degenerate.
\bigskip
\proclaim Theorem 4.2.
Let $\L\iN{H^*}$ be such that
$\o^*(\L)=\L$.  There exists an element $\sgma\iN H^0$ such that 
the twining character of the Verma 
module of highest weight $\L$ with respect to $\tilde\o$ is
  $$\N_{\L}^{\tilde\o}(h)=e(\L(h_{\tom}))\,\N_{\L}^{\o}(h-h_{\tom})\,,\ \
  h\in H\,.$$
If moreover, $\L$ is also integrable, then the twining character of the 
irreducible module of highest weight $\L$ with respect to $\tilde\o$ is
  $$\Psi_{\L}^{\tilde\o}(h)=e(\L(h_{\tom}))\,\Psi_{\L}^{\o}(h-h_{\tom})\,,\ \
  h\in H\,.$$
\par\noindent
{\bf Proof.} We write the irreducible twining characters for 
the diagram automorphism $\o$ as
  $$  \Psi_\Lambda^\omega = \sum_{\lambda\leq\Lambda} m_\lambda^{\omega}\,
  e(\lambda)\eqno(6) $$
and for the generalized diagram automorphism as
  $$  \Psi_\Lambda^{\tom} = \sum_{\lambda\leq\Lambda} m_\lambda^{\tom}\,
  e(\lambda) \, . \eqno(7) $$
Let $\l$ be an element in $H^*$ such that
$m_{\l}^{\o}\not=0$ (or equivalently, $m_\lambda^{\tom}\neq 0$). Then
\hfill\break
$\l=\L-\sum_{i\in I} k_i\a_i$, where for each $i\iN I$, $k_i$ is a 
non-negative integer, $k_{\do i}=k_i$, and
$k_i=0$ unless $i$ is $\do$-conjugate to an element in $\cI$. We find that
  $$ m_\lambda^{\tom} =  m_\lambda^{\omega} \prod_{i\in I} (\xi_i)^{k_i} 
  = m_\lambda^{\omega} \prod_{i\in\cI} \llb\prod_{l=0}^{N_i-1} \xi_{\do^l i}
  \lrb^{k_i} \, . $$
For each $i\iN I$, we define $\sigma_i\iN\Bbb C$ to be as follows:
\be \prod_{l=0}^{N_i-1} \xi_{\do^l i} = \eE^{\sigma_i/N_i} \,. \be
The imaginary part of $\sigma_i$ is of course only fixed modulo 
$2\pi N_i$, and we may put $\sigma_{\do i}=\sigma_i$ for all $i\iN I$.  
This allows us to express $m_\lambda^{\tom}$ as
  $$ m_\lambda^{\tom} = m_\lambda^{\omega} \prod_{i\in\cI} 
  \eE^{\sigma_ik_i/N_i} 
  = m_\lambda^{\omega} \prod_{i\in I} \eE^{\sigma_i k_i}  \, . $$

For $i\iN I$, let $\L_i$ denote the fundamental weights, satisfying
$(\L_i,\a_j)=\d_{ij}$ for all $i,j\iN I$. We can write $k_i$ as
$k_i = (\Lambda-\lambda, \Lambda_{i})\,.$ Define the element
  $$\sigma:= \sum_i \sigma_i \Lambda_{i} \,.$$
Then $\sigma$ is in $(H^*)^0$. We obtain 
  $$m_\lambda^{\tom} = m_{\l}^{\o}\, \eE^{(\Lambda-\lambda,\sigma)}\, . $$
Let $\varphi\!:\,H^*\to H$ be the bijection induced by the bilinear form
$(.\,,.)$ on $H$, and let $\sgma$ denote $\varphi(\sigma)$. 
Substituting the above expression in $(7)$ and using $(6)$ we find
that $\Psi_\Lambda^{\tom}(h) = e(\Lambda(\sgma))\Psi_\Lambda^\omega(h-
\sgma)\,,$
proving the Theorem.
\hfill\pf
\bigskip
Hence the effect of a generalized diagram automorphism,
as compared to the corresponding ordinary diagram automorphism, consists in
a shift in the argument
and a multiplication by an overall factor. In case $\tom$ has finite order,
this factor is of course a phase. Note that the imaginary part of $\sigma$ 
is defined only up to $2\pi$ times an element of the 
weight lattice of $G$; the real part of $\sigma$ is unique, however; it is zero
if $\tom$ has finite order.
\bigskip\bigskip
\def\comp  {{\it Commun.\ Math.\ Phys.}}
\def\joal  {{\it Journal of Algebra}}
\def\nupb  {{\it Nuclear Physics B}}
\def\ijmp  {{\it Intern.\ J.\ of Modern Physics A}}
\def\newl  {\hfill\break\noindent\hbox{\phantom{xxi}}}
\centerline{\bf References}
\bigskip\noindent
1. R.E.\ Borcherds, Generalized \kma s. \joal\ {\bf 115} (2), (1988)
\newl 501--512 
\bigskip\noindent
2. J.\ Fuchs, A.N.\ Schellekens, and C.\ Schweigert, 
{}From Dynkin diagram symmetries to\newl fixed point structures.
{\it Preprint hep-th/9506135}, \comp, in press
\bigskip\noindent
3. J.\ Fuchs, A.N.\ Schellekens, and C.\ Schweigert, 
The resolution of field identification fixed\newl points in diagonal coset 
theories. \nupb\ {\bf 461} (1996) 371--406
\bigskip\noindent
4. J.\ Fuchs, A.N.\ Schellekens, and C.\ Schweigert, 
A matrix $S$ for all simple current exten-\newl sions.
{\it Preprint hep-th/9601078}, \nupb, in press
\bigskip\noindent
5. V.G.\ Kac, ``Infinite dimensional Lie Algebras'', third ed.,
Cambridge University Press\newl 1990
\bigskip\noindent
6. A.N.\ Schellekens and S.\ Yankielowicz, 
Simple currents, modular invariants, and fixed\newl points,
\ijmp\ {\bf 5} {(1990)} 2903--2952           

\end